\documentclass[11pt]{article}

\usepackage{graphicx,enumitem}
\usepackage{graphicx,psfrag,epsf}
{%
\centering\addtocounter{figure}{1}
\begin{enumerate}[%
itemsep=2pt,parsep=0em,
label={(\alph*)},
ref={\thefigure.(\alph*)}
]}%
{\end{enumerate}\addtocounter{figure}{-1}}

\usepackage[margin=1in]{geometry}

\usepackage[utf8]{inputenc}
\usepackage{amsmath}
\usepackage[english]{babel}
\usepackage{amssymb}
\usepackage{amsfonts}
\usepackage{bm}
\usepackage{graphicx}
\usepackage{amsmath}
\usepackage{amsthm}
\usepackage{xcolor}
\usepackage[hidelinks=true]{hyperref}
\usepackage{parskip}
\usepackage{multicol}
\usepackage{float}

\usepackage{verbatim}

\floatstyle{plaintop}
\restylefloat{table}
\usepackage[tableposition=top]{caption}

\usepackage{rotating}
\usepackage{subcaption}
\usepackage{multirow}

\usepackage{chngcntr}
\usepackage{apptools}
\AtAppendix{\counterwithin{lem}{subsection}}
\AtAppendix{\counterwithin{proposition}{subsection}}
\AtAppendix{\counterwithin{equation}{section}}

\newtheoremstyle{mystyle}%
  {7pt}   
  {7pt}   
  {\itshape}  
  {}       
  {\bfseries} 
  {.}      
  {0.5em}  
  {}

\theoremstyle{mystyle}

\newtheorem{thm}{Theorem}[section]

\newtheorem{proposition}{Proposition}[section]

\newtheorem{ejemplo}{Example}[section]

\usepackage{booktabs}
\usepackage{animate}

\usepackage{natbib}
\hypersetup{urlcolor=blue, citecolor=red, colorlinks=true, linkcolor=red} 

\numberwithin{equation}{section}
\numberwithin{figure}{section}
\numberwithin{table}{section}

\usepackage{authblk}

\usepackage{blindtext}
\usepackage{lineno}

\begin{document}

\def\spacingset#1{\renewcommand{\baselinestretch}%
{#1}\small\normalsize} \spacingset{1}


\title{\bf Non-Stationary Covariance Functions for Spatial\\ Data on Linear Networks}

\author{
Alfredo Alegr\'ia\thanks{Email: alfredo.alegria@uc.cl}
}

\date{}

\maketitle

\vspace{-1.5em}
\begin{center}
{\small
Departamento de Estad\'istica, Pontificia Universidad Cat\'olica de Chile,\\ Santiago, Chile
}

\vspace{2em}
{\small\itshape \today}
\vspace{2.5em}
\end{center}

\begin{abstract}
\noindent We introduce a novel class of non-stationary covariance functions for random fields on linear networks that allows both the variance and the correlation range of the random field to vary spatially. The proposed covariance functions are useful to model random fields with a spatial dependence that is locally isotropic with respect to the resistance metric, a distance that reflects the topology of the network. The framework admits explicit stochastic representations of the associated random fields and can be naturally extended to matrix-valued covariance functions for vector-valued random fields. We assess the statistical and computational performance of a weighted local likelihood estimator for the proposed models using synthetic data generated on the street network of the University of Chicago neighborhood. \\

\end{abstract}

\section{Introduction}

Data indexed by locations on a network arise naturally in a variety of applications, such as processes on road, river, and electrical networks. Recent developments in the theory and methodology of spatial statistics on linear networks can be found in \cite{baddeley2017stationary}, \cite{anderes2020isotropic}, \cite{cronie2020inhomogeneous}, \cite{Moradi}, \cite{baddeley2021analysing}, \cite{porcu2023stationary}, \cite{bolin2024gaussian}, \cite{moller2022cox}, \cite{moradi2024summary},   \cite{filosi2025vector},  \cite{gonzalez2025nonparametric}, \cite{bolin2025new}, and \cite{Alegría20102025}. 

Specifically, a linear network $\mathcal{L}$ is the union $\mathcal{L} = \bigcup_{i=1}^m \ell_i$ of closed line segments $\ell_1,\hdots,\ell_m$ in $\mathbb{R}^2$, with $m<\infty$. Each segment has positive finite length and intersects other segments only at endpoints. Throughout, the network is assumed to be path-connected. 

A second-order Gaussian random field (GRF), namely $\{X(s): s \in \mathcal{L}\}$, is completely characterized by its mean function, which we assume to be zero for simplicity, and by its covariance function, given by 
$$K(s,t) = {\rm cov}[ X(s), X(t)], \qquad s,t \in \mathcal{L}.$$ 
The construction of valid and versatile covariance functions is challenging because of the requirement that the covariance matrix of $(X(s_1),\hdots,X(s_N))^{\top}$, given by
$$  \bm{\Sigma} = \big[ K(s_i,s_j) \big]_{i,j=1}^N,    $$
must be positive semidefinite, for all $N \in \mathbb{N}$ and for all collection of spatial locations $s_1, \dots, s_N$ in $\mathcal{L}$. In other words, $K$ must be a positive semidefinite function. The spectral theory of random fields on Euclidean spaces is useful to check for positive semidefiniteness of a candidate function (see, e.g., \citealp{yaglom1987correlation}). In the linear network setting, however, analogous characterizations are unavailable.  

It is commonly assumed that the covariance function is isotropic, meaning that there exists a function $C:[0,\infty) \rightarrow \mathbb{R}$ such that
 $$K(s,t)=C(\text{d}(s,t)), \qquad s,t\in\mathcal{L},$$  with $\text{d}(\cdot,\cdot)$ standing for a distance metric on $\mathcal{L}$. In practice, treating network coordinates as ordinary Euclidean coordinates can distort the proximity between locations and produce unrealistic results. To address this, \cite{anderes2020isotropic} established positive definiteness for certain classes of isotropic covariance functions under the resistance metric (an extension of the resistance metric from electrical network theory) and showed that, under the geodesic (shortest-path) metric, these classes are valid only under restrictive assumptions on the network structure. Both the resistance and geodesic metrics, unlike the Euclidean metric, take the network topology into account.

Although isotropy is a convenient assumption, it can be overly restrictive. In geostatistical applications, there is increasing demand for more flexible models, particularly when analyzing data over large domains where the properties of the process may vary across space. See \cite{bolin2025new} for a recent development of non-stationary random fields on metric graphs formulated as solutions to stochastic partial differential equations, and an application to traffic speed data. In the context of log Gaussian Cox processes (LGCPs) \citep{cox1955some}, where the intensity of point occurrences is governed by a log Gaussian random field, isotropy of the underlying GRF prevents the number of points and the spatial extent of point clusters from varying across the domain. This limitation has been studied for point processes in Euclidean spaces by \cite{dvovrak2019quick}. Similar issues arise on linear networks, as discussed by \cite{baddeley2017stationary}, encouraging the development of more flexible models.

Motivated by this, we propose a general class of non-stationary covariance functions on linear networks, that is, covariance functions not solely determined by distance. The proposed covariance models allow both the variance and the correlation range of the random field to vary across the network, providing a network analogue of models previously proposed in Euclidean spaces \citep{paciorek2003nonstationary,paciorek2006spatial,anderes2011local,kleiber2012nonstationary,mateu2013class,liang2016class}. Our framework produces covariance functions that are locally isotropic with respect to the resistance metric, ensuring interpretability and alignment with the topology of the network. 

We exploit the connection between the proposed covariance functions and the class of generalized Stieltjes functions (see \citealp{koumandos2019generalized,menegatto2020positive} for a detailed treatment) to provide examples of non-stationary covariance functions in closed form.

The proposed covariance functions admit explicit stochastic representations of the associated random fields. We discuss the importance of these representations in the design of computationally efficient geostatistical simulation algorithms. In addition, the framework is extended to construct matrix-valued covariance functions for vector-valued GRFs on linear networks, a relevant setting for analyzing multiple georeferenced variables simultaneously.

To address statistical inference, we revisit the weighted local likelihood approach, originally introduced by \cite{anderes2011local} for Euclidean spaces, and adapt it to the network setting. We then conduct a simulation study using synthetic data generated on the street network of the University of Chicago neighborhood, in which we estimate the local correlation range, examine the behavior of the method under different weighting schemes, and evaluate its computational performance.

The paper is organized as follows. Section \ref{sec:preliminaries} reviews the definition and properties of the resistance metric. Section \ref{sec:results} presents the main theorem, examples, and illustrations. In Section \ref{sec:developments} we provide a stochastic representation and extend the framework to construct matrix-valued covariance functions. Section \ref{sec:simulation} describes the weighted local likelihood method and reports the results of the simulation study. Section \ref{sec:discussion} concludes the paper with a discussion and directions for future research. Proofs are provided in the appendices.

\section{Resistance metric}
\label{sec:preliminaries}

This section reviews the definition and properties of the resistance metric introduced in \cite{anderes2020isotropic}, as it is essential for building the non-stationary covariance functions proposed in this work. 

We begin by representing the linear network as a metric graph $(\mathcal{V}, \mathcal{E})$, where the vertex set $\mathcal{V} = \{v_1, \dots, v_n\}$ corresponds to the endpoints of $\mathcal{L}$, and the edge set $\mathcal{E}$ corresponds to the interiors of the linear segments $\ell_1, \dots, \ell_m$. The length of $\ell_i$ is denoted by $|\ell_i|$. We write $v_i \sim v_j$ if the distinct vertices $v_i$ and $v_j$ are connected by an edge. For any pair of vertices $v_i, v_j \in \mathcal{V}$, the conductance function is defined as 
$$
\text{con}(v_i,v_j) =
\begin{cases}
    1/\|v_i - v_j\|, & \text{if } v_i \sim v_j,\\
    0, & \text{otherwise.}
\end{cases}
$$  
Using this function, we can build a matrix $\bm{\Delta}$ indexed by the vertices in $\mathcal{V}$, with entries  
$$
\bm{\Delta}(v_i,v_j) =
\begin{cases}
    1 + c(v^*), & \text{if } v_i = v_j = v^*,\\
    c(v_i), & \text{if } v_i = v_j \neq v^*,\\
    -\text{con}(v_i,v_j), & \text{otherwise,}
\end{cases}
$$  
where $v^*\in\mathcal{V}$ is an arbitrary reference vertex and $c(v_i) = \sum_{v \sim v_i} \text{con}(v_i,v)$. The resulting matrix can be viewed as a modification of the graph Laplacian, where the inclusion of an additional unit term on the diagonal guarantees strict positive definiteness.  This construction is invariant to the choice of $v^*$ \citep{anderes2020isotropic}.

To define the resistance metric, we introduce an auxiliary GRF on $\mathcal{L}$, given by
\begin{equation}
\label{auxiliary}
    A(s) =  A_0(s) + \sum_{i=1}^m A_i(s), \qquad s\in \mathcal{L},
\end{equation}
where $A_0,\hdots,A_m$ are independent GRFs defined as follows:
\begin{itemize}
    \item The base component $A_0$ is first defined on the vertices as a multivariate normal vector: 
$$(A_0(v_1),\hdots,A_0(v_n))^\top \sim \mathcal{N}_n\big(\bm{0}, \bm{\Delta}^{-1}\big).$$
 Next, $A_0$ is extended from the vertices to the edges using linear interpolation.

    \item For $i=1,\hdots,m$, 
$$
A_i(s) =
\begin{cases}
    B_i(s), & \text{for } s \in \ell_i,\\
    0, & \text{for } s \in \mathcal{L} \setminus \ell_i,
\end{cases}
$$
where $B_1, \dots, B_m$ are independent standard Brownian bridges, each defined on an interval of length $|\ell_i|$, starting and ending at zero, and independent of $A_0$. Since each $B_i$ is defined on an interval, writing $B_i(s)$ implicitly identifies points $s$ along $\ell_i$ with points in this interval via arc length parameterization. 
\end{itemize}

For additional details on the construction of the auxiliary GRF, we refer the reader to \cite{anderes2020isotropic} and \cite{moller2022cox}. 

The resistance metric on $\mathcal{L}$ can now be defined as the variogram of the auxiliary GRF in Equation (\ref{auxiliary}), that is,
\begin{equation}
    \label{metric}
    \text{d}_R(s,t) = {\rm var}[A(s) - A(t)], \qquad s,t\in \mathcal{L}.
\end{equation}
This metric extends the classical resistance metric from electrical network theory \citep{klein}, originally defined on the vertex set, to the continuum of points along the edges of the graph. When evaluated on any additional points along the edges, $\text{d}_R$ coincides with the metric obtained from the corresponding discrete electrical network. Interestingly, the resistance and geodesic distances coincide when $\mathcal{L}$ is a tree, that is, a connected graph with no cycles. \cite{moller2022cox} provide more fully developed algebraic expressions for the resistance metric, which facilitate its implementation.

We turn to a description of the range of isotropic covariance functions attainable using the resistance metric. A function $\varphi : [0,\infty) \to \mathbb{R}$ is said to be {completely monotone} on $[0,\infty)$ if it is continuous on $[0,\infty)$, infinitely differentiable on $(0,\infty)$, and satisfies $(-1)^j \varphi^{(j)}(z) \ge 0$, $z > 0$, for every integer $j \ge 0$, where $\varphi^{(j)}$ denotes the $j$-th derivative of $\varphi$ and $\varphi^{(0)} = \varphi$. \cite{anderes2020isotropic} showed that for any completely monotone function $\varphi$, the function $K(s,t)=\varphi\big(\text{d}_R(s,t)\big)$ is (strictly) positive definite on $\mathcal{L} \times \mathcal{L}$. As a consequence, many popular parametric classes of covariance functions developed for Euclidean spaces are also valid with respect to the resistance metric for general linear networks (see \citealp[Table 1]{anderes2020isotropic}). 

In contrast, the geodesic distance provides a more restrictive alternative. For example, if $\mathcal{L}$ contains three distinct paths connecting two points $s,t \in \mathcal{L}$, the classical exponential covariance function, defined with respect to the geodesic metric, is no longer positive definite on this network \citep{anderes2020isotropic}. Consequently, covariance functions defined in terms of completely monotone functions, which are mixtures of exponential covariance functions, are no longer positive definite on such networks when the geodesic metric is used. This considerably restricts the use of the geodesic distance, as many real situations, such as phenomena on street networks, typically do not satisfy this condition. 

\cite{moller2022cox} also emphasize, in the context of point processes on linear networks, the advantages of the resistance metric over the geodesic metric in terms of flexibility. In particular,  the use of $\text{d}_R$ allows one to model stronger clustering effects in Cox processes and higher degree of repulsiveness in determinantal point processes.

\section{Non-stationary covariance functions}
\label{sec:results}

\subsection{Main result}

We now present a general theorem for constructing non-stationary families of covariance functions on linear networks.

\begin{thm}
\label{teo1}
Let $a: \mathcal{L} \rightarrow (0,\infty)$ and $b: \mathcal{L} \rightarrow (0,\infty)$, and define $\alpha(s,t)= [a(s) + a(t)]/2$.
 Let $\text{d}_R(s,t)$ be the resistance metric between $s$ and $t$, and let $F$ denote a positive finite measure on $(0,\infty)$. Then, the function 
\begin{equation}
    \label{cova}
    K(s,t) =   \displaystyle \int_0^\infty \frac{\sqrt{b(s)b(t)}}{\sqrt{ \alpha(s,t) + \text{d}_R(s,t)w}} \text{d}F(w), \qquad s,t\in\mathcal{L},
\end{equation}
is positive semidefinite on $\mathcal{L}\times \mathcal{L}$.
\end{thm}

Our proof of Theorem \ref{teo1}, presented in Appendix \ref{app1}, adapts ideas introduced in \cite{paciorek2003nonstationary} for establishing positive definiteness of non-stationary covariance functions in Euclidean spaces. To the best of our knowledge, the present work is the first to extend this approach to linear networks. The fact that the resistance metric is defined as the variogram of a specific auxiliary GRF gives the proof a distinctive character.

 Observe that if $s$ and $t$ are near a fixed location $s_0\in \mathcal{L}$, and both functions $b$ and $a$ are sufficiently smooth, then 
\begin{equation*}
\label{cova_aprox}
K(s,t) \approx    \int_0^{\infty} \frac{b(s_0)}{\sqrt{a(s_0) + \text{d}_R(s,t) w}}  \text{d}F(w). 
\end{equation*}
 In other words, $K(s,t)$ depends approximately only on the resistance metric $\text{d}_R(s,t)$. This property allows us to construct covariance models with locally isotropic behavior with respect to $\text{d}_R$, reflecting the topology of the linear network. 
 
 In this framework, the function $b$ provides a way to modulate the locally varying variance of the random field, while the function $a$ regulates how quickly the covariance function decreases with distance. In particular, $a$ can be adjusted to influence the locally adaptive correlation range (the distance beyond which correlations fall below $0.05$).

\subsection{Examples}
\label{sec:examples}

A direct calculation shows that the covariance function in Equation (\ref{cova}) admits the following representation:
$$  K(s,t) =  \begin{cases}
\displaystyle \frac{b(s)}{\sqrt{a(s)}}  F((0,\infty)), & \text{ if } s=t, \\   \\ 
\displaystyle   \sqrt{\frac{b(s)b(t)}{\text{d}_R(s,t)}}  \psi\left(  \frac{\alpha(s,t)}{\text{d}_R(s,t)}  \right),  & \text{ if } s\neq t,
 \end{cases}
 $$
where 
\begin{equation}
\label{stieltjes}
\psi(r) = \int_0^\infty \frac{1}{\sqrt{r+w}}\text{d}F(w), \qquad r>0. 
\end{equation}
Functions of the form (\ref{stieltjes}) belong to the class of generalized Stieltjes functions of order $1/2$ (see, e.g., \citealp{koumandos2019generalized}). When $a$ and $b$ are constant functions, we recover a subclass of the isotropic covariance functions presented in \cite{anderes2020isotropic}, since every generalized Stieltjes function is completely monotone (see, e.g., \citealp[Section 2.2]{menegatto2020positive}).

Since closed-form expressions for generalized Stieltjes functions are available in standard tables of integral transforms such as \cite[Chapter 14]{erdelyi}, it is straightforward to derive explicit examples of non-stationary covariance functions.

\begin{ejemplo}
Let $F$ be the Dirac delta measure at $\eta>0$, that is, $\text{d}F(w) = \text{d}\delta_{\eta}(w)$. Then, we immediately obtain the covariance function
\begin{equation}
    \label{cauchy}
    K(s,t)  = \sqrt{\frac{b(s) b(t)}{\alpha(s,t)}}   \left(  1  + \frac{\eta \text{d}_R(s,t)}{\alpha(s,t)} \right)^{-1/2}.
\end{equation}
\end{ejemplo}

\begin{ejemplo}
Let $\text{d}F(w) =  (1/\eta) \exp(-w/\eta) \text{d}w$, for $\eta>0$. Then, using Equation 14.4.10 in \cite{erdelyi}, we obtain
\begin{equation*}
    \label{erfc}
    K(s,t)  = \begin{cases} 
   \displaystyle  \frac{b(s)}{\sqrt{a(s)}},   &  \text{ if }  s=t,   \\ \\ 
   \displaystyle  \sqrt{\frac{\pi b(s)b(t)}{\eta}} \exp\left(\frac{\alpha(s,t)}{\eta \text{d}_R(s,t)}\right) \text{erfc}\left(\sqrt{\frac{\alpha(s,t)}{\eta \text{d}_R(s,t)}} \right),  &   \text{ if } s\neq t,
    \end{cases}
\end{equation*}
where erfc denotes the complementary error function. 
\end{ejemplo}

\begin{ejemplo}
Let $\text{d}F(w) = \sqrt{\eta} \pi^{-1} w^{-1/2}(\eta+w)^{-1}  \text{d}w$, for $\eta>0$. This is a scaled Beta prime probability measure. Then, using Equation 14.4.9 in \cite{erdelyi}, we obtain
\begin{equation*}
    \label{hyperg}
    K(s,t)  = \begin{cases} 
   \displaystyle  \frac{b(s)}{\sqrt{a(s)}},   &   \text{ if }  s=t,   \\ \\
   \displaystyle \frac{2}{\pi} \sqrt{\frac{b(s)b(t)}{\eta}} \, _{2}F_1\left(1,1/2,3/2; 1 - \frac{\alpha(s,t)}{\eta \text{d}_R(s,t)}\right),  &   \text{ if }  s\neq t,
    \end{cases}
\end{equation*}
where $_{2}F_1$ denotes the Gauss hypergeometric function. 
\end{ejemplo}

 \subsection{Illustration}
\label{sec:illustration}

 To illustrate our proposal, we simulate a GRF on a linear network formed by the mortar lines of a brick wall. This network, consisting of 156 vertices and 203 edges, and available in the \texttt{R} package \texttt{spatstat} \citep{baddeley2005spatstat}, has been the focus of studies on the positions of spider webs of the urban wall spider Oecobius navus \citep{voss2007habitat}. 
 
\begin{figure}
    \centering
    \includegraphics[scale=0.8]{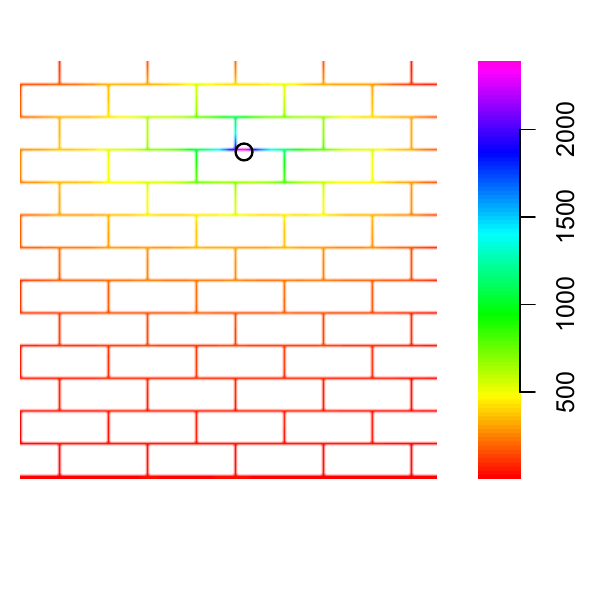}
        \includegraphics[scale=0.8]{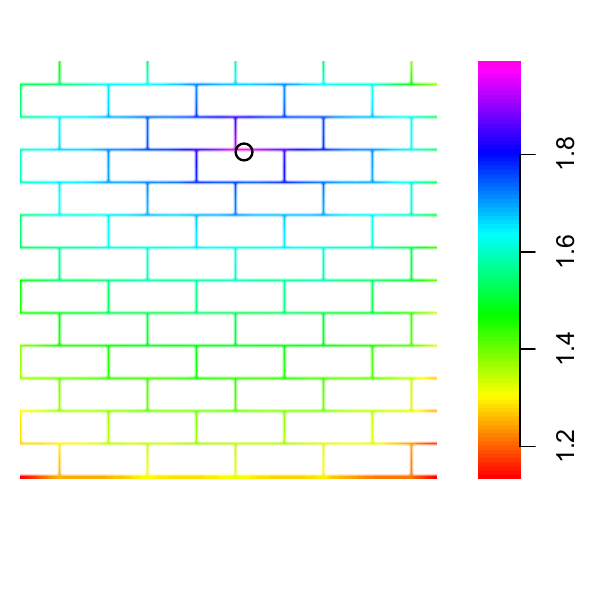}
    \caption{Functions $a$ (left) and $c$ (right), defined in Equations (\ref{eq:aa}) and (\ref{eq:bb}), respectively, on the \texttt{spiders} linear network. The circle indicates the point $s^{*}\in\mathcal{L}$ involved in the definitions of $a$ and $c$.}
    \label{fig:functions}
\end{figure}

 We use the covariance function in Equation (\ref{cauchy}) with $\eta=1$ and the parameterization $b(s)=c(s)\sqrt{a(s)}$, for some function $c:\mathcal{L}\rightarrow (0,\infty)$, which ensures that ${\rm var}[X(s)] = K(s,s) = c(s)$, for all $s\in \mathcal{L}$. Specifically, we consider
\begin{equation}
\label{eq:aa}
a(s) =  2500 \times \exp\left( - \frac{\text{d}_R(s,s^{*})}{100} \right), \qquad s\in\mathcal{L},
\end{equation}
 and 
\begin{equation}
\label{eq:bb}
 c(s) = 2 \times \exp\left( - \frac{\text{d}_R(s,s^{*})}{1000} \right), \qquad s\in\mathcal{L},
 \end{equation}
 where $s^{*}$ is a fixed point on $\mathcal{L}$. Both functions attain their maximum at $s^{*}$ and decay exponentially as $s$ moves away from this point, as illustrated in Figure \ref{fig:functions}.  Accordingly, we expect the simulated GRF to exhibit larger correlation range and higher variance at locations near $s^{*}$, decreasing as one moves away from this point. 
 
\begin{figure}
    \centering
    \includegraphics[scale=0.75]{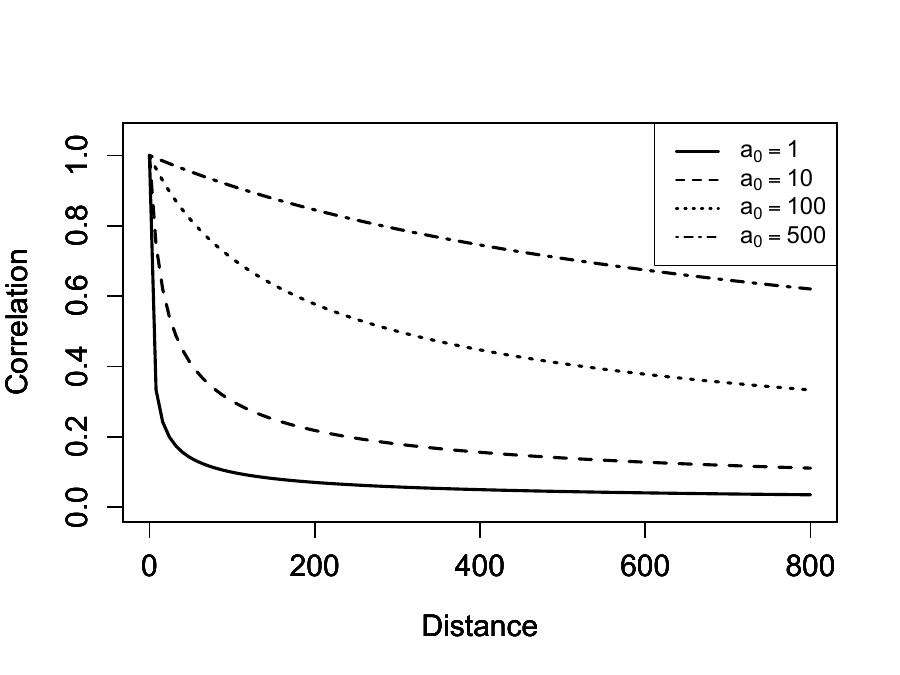}
    \caption{Mapping $\text{d}_R \mapsto \sqrt{a_0}/\sqrt{a_0 + \text{d}_R }$ for selected values of $a_0$, illustrating how the local isotropic correlation structure of the non-stationary covariance function (\ref{cauchy}) depends on $a_0$. }
    \label{fig:corr}
\end{figure}

 To better illustrate the local behavior of this covariance function, Figure \ref{fig:corr} shows the curves 
 $$\text{d}_R \mapsto   \frac{\sqrt{a_0}}{\sqrt{a_0 + \text{d}_R }}$$
 for different values of $a_0$. These curves capture part of the local isotropic correlation structures that the  non-stationary covariance function exhibits. They are displayed over the interval $(0,800)$ since the maximum resistance distance between two points in this domain is approximately 837.
 
\begin{figure}
    \centering
    \includegraphics[scale=0.72]{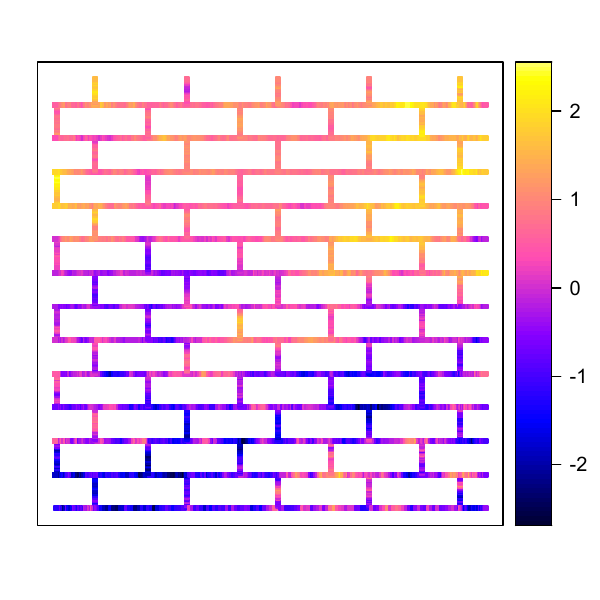}
        \includegraphics[scale=0.88]{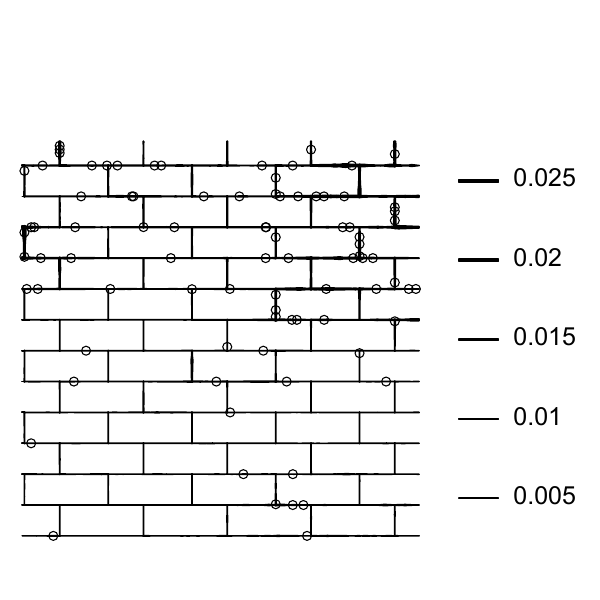}
    \caption{(Left) Simulated GRF on the \texttt{spiders} linear network, with the non-stationary covariance function given in Equation (\ref{cauchy}). (Right)  Corresponding realization of a LGCP, with line thickness proportional to the magnitude of the log Gaussian field.}
    \label{fig:illustration}
\end{figure}

Figure \ref{fig:illustration} (left panel) shows the simulated GRF over 5,684 sites (28 per edge). The simulated field attains high values at sites in the top part of the linear network, and because the correlation range is greater in this region, these high values extend over a wide area. In contrast, at sites in the lower part of the linear network, the simulated field varies more rapidly, reflecting a smaller correlation range. 

The right panel of Figure \ref{fig:illustration} shows the realization of the corresponding LGCP, where the log Gaussian random field is defined as $\Lambda(s) = 0.002 \times \exp\{X(s)\}$, with $X(s)$ denoting the underlying GRF. Recall that, conditioned on $\Lambda$, this point process on $\mathcal{L}$ is almost surely a Poisson process with intensity function $\Lambda$ (with respect to the arc length measure on $\mathcal{L}$). As expected, more points occur in the upper part of the network, covering a larger area, while fewer points are present in the lower part. For a more detailed treatment of Cox processes on linear networks, we refer the reader to \cite{moller2022cox}.

\section{Extensions of the framework}
\label{sec:developments}

\subsection{Stochastic representation}

The following proposition provides an explicit construction of an elementary random field whose covariance function is given by (\ref{cova}). Throughout this subsection, we assume without loss of generality that $F$ is a probability measure on $(0,\infty)$.

\begin{proposition}
\label{prop:stochastic}
Let $Z$ be a real-valued random variable with probability density function $g$. Let $W$ be a positive random variable with distribution $F$. Let $\epsilon$ be a random variable with zero mean and unit variance. Let $\{A(s): s\in\mathcal{L}\}$ denote the auxiliary GRF defined in Equation (\ref{auxiliary}). Assume that $Z$, $W$, $\epsilon$, and $A$ are mutually independent. Let $\phi(\cdot; \mu, \sigma^2)$ denote the probability density function of a normal distribution with mean $\mu$ and variance $\sigma^2$. Then, the random field
\begin{equation}
\label{stochastic}
X(s) =  \epsilon \, \sqrt{\frac{b(s)}{g(Z)}} \left( \frac{\pi}{W} \right)^{1/4}  \phi\left( Z; \frac{A(s)}{\sqrt{2}}, \frac{a(s)}{4W} \right), \qquad s\in\mathcal{L},  
\end{equation}
has zero mean and covariance function given by (\ref{cova}).
\end{proposition}

A proof of this result is provided in Appendix \ref{app2}.

The stochastic representation introduced in Equation (\ref{stochastic}) extends Proposition 3.4 of \cite{Alegría20102025} from the isotropic setting to the non-stationary case. 

The importance of the stochastic representation (\ref{stochastic}) lies in its potential for developing computationally efficient algorithms to simulate GRFs with covariance functions of the form (\ref{cova}). Recall that simulation based on the Cholesky decomposition of the covariance matrix has cubic computational complexity, which motivates the development of more efficient methods (see \citealp{lantuejoul2001geostatistical} for a comprehensive treatment in Euclidean spaces).

 Since the random field in (\ref{prop:stochastic}) is not Gaussian, one strategy is to add several independent copies of this elementary random field, rescale appropriately,  and then invoke a central limit approximation. More precisely, an approximately GRF with covariance function (\ref{cova}) is given by
 $$ X(s) = \frac{1}{\sqrt{J}}\sum_{j=1}^J   \epsilon_j  \, \sqrt{\frac{b(s)}{g(Z_j)}} \left( \frac{\pi}{W_j} \right)^{1/4}  \phi\left( Z_j; \frac{A^{(j)}(s)}{\sqrt{2}}, \frac{a(s)}{4W_j} \right), \qquad s\in\mathcal{L},  
  $$
where $J$ is a large positive integer, and $\epsilon_1,\hdots,\epsilon_J$, $Z_1,\hdots,Z_J$,  $W_1,\hdots,W_J$, and  $A^{(1)},\hdots,A^{(J)}$ are independent copies of the random ingredients involved in Proposition \ref{prop:stochastic}.

We refer the reader to \cite[Section 4.5]{Alegría20102025} for a more detailed treatment of the computational aspects of these types of algorithms in the linear network scenario, as well as for implementation guidelines.

\subsection{Vector-valued GRFs}

In multivariate spatial statistics, several georeferenced variables are analyzed simultaneously \citep{wackernagel2003multivariate}.
Interest focuses not only on the marginal covariance functions of each variable but also on the cross-covariance functions between them. 

 Formally, consider a vector-valued GRF, namely $\{\bm{X}(s) = (X_1(s),\hdots,X_q(s))^\top: s\in\mathcal{L}\}$, where $q$ is the number of components, and ${\rm var}[X_i(s)] < \infty$, for all $s\in\mathcal{L}$ and $i=1,\hdots,q$. In this setting, the covariance function is matrix-valued, $\bm{K}:\mathcal{L}\times \mathcal{L}\rightarrow \mathbb{R}^{q \times q}$, with entries 
$$  K_{ij}(s,t) =  {\rm cov}\big[X_i(s),X_j(t)\big], \qquad s,t\in\mathcal{L}, \,\,  i,j=1,\hdots,q.$$ 
The covariance function must be positive semidefinite in the sense that, for each $N\in\mathbb{N}$ and any collection of locations $s_1,\hdots,s_N \in\mathcal{L}$, the covariance matrix of the $Nq$-dimensional vector $(\bm{X}(s_1)^\top, \hdots, \bm{X}(s_N)^\top)^\top$ is positive semidefinite.

The following theorem provides a construction of non-stationary matrix-valued covariance functions, extending Theorem \ref{teo1} to the vector-valued setting.
\begin{thm}
\label{teo2}
Let $a_i,b_i: \mathcal{L} \rightarrow (0,\infty)$, for $i=1,\hdots,q$, and define $\alpha_{ij}(s,t)= [a_i(s) + a_j(t)]/2$. Let $\bm{\beta} = [\beta_{ij}]_{i,j=1}^q$ be a symmetric positive semidefinite matrix. Let $\text{d}_R(s,t)$ denote the resistance metric between $s$ and $t$, and let $F$ be a positive finite measure on $(0,\infty)$. Then, the matrix-valued function $\bm{K}: \mathcal{L}\times\mathcal{L}\rightarrow \mathbb{R}^{q\times q}$, with entries 
\begin{equation}
    \label{cova_mult}
    K_{ij}(s,t) =   \displaystyle \int_0^\infty \frac{\beta_{ij} \sqrt{b_i(s)b_j(t)}}{\sqrt{ \alpha_{ij}(s,t) + \text{d}_R(s,t)w}} \text{d}F(w), \qquad s,t\in\mathcal{L}, \,\,  i,j=1,\hdots,q,
\end{equation}
is positive semidefinite.
\end{thm} 
In Appendix \ref{app3}, we provide a proof of Theorem \ref{teo2}. This theorem is the network counterpart of a result presented in \cite{kleiber2012nonstationary} for Euclidean spaces. 

In this construction, the functions $b_i$ and $a_i$ allow us to control, respectively, the local variance and correlation range of the $i$-th component of the vector-valued GRF, in a similar fashion to the univariate ($q=1$) case. The parameters $\beta_{ij}$, in turn, allow us to modulate the colocated correlation coefficient between the $i$-th and $j$-th components. More precisely, the effect of $\beta_{ij}$ on the colocated correlation is given by
$$ {\rm cor}[X_i(s),X_j(s)] = \frac{{\rm cov}[X_i(s),X_j(s)]}{\sqrt{ {\rm var}[X_i(s)]  {\rm var}[X_j(s)]  }} = \frac{\beta_{ij}}{\sqrt{\beta_{ii}\beta_{jj}}}  \frac{  \left[  a_i(s)a_j(s) \right]^{1/4}  }{  \sqrt{\alpha_{ij}(s,s)}  }, $$
for each $s\in\mathcal{L}$ and $i\neq j$.

Using the representation in terms of generalized Stieltjes functions of order $1/2$, as described in Section \ref{sec:examples},  the function (\ref{cova_mult}) can be expressed as
$$  K_{ij}(s,t) =  \begin{cases}
\displaystyle \beta_{ij} \sqrt{\frac{b_i(s)b_j(s)}{\alpha_{ij}(s,s)}}  F((0,\infty)), & \text{ if } s=t, \\   \\ 
\displaystyle \beta_{ij}  \sqrt{\frac{b_i(s)b_j(t)}{\text{d}_R(s,t)}}  \psi\left(  \frac{\alpha_{ij}(s,t)}{\text{d}_R(s,t)}  \right),  & \text{ if } s\neq t,
 \end{cases}
 $$
with $\psi$ defined according to Equation (\ref{stieltjes}). Thus, analogous examples to those presented in Section \ref{sec:examples} can be constructed for the vector-valued case.

In particular, when $a_i$ and $b_i$ are constant functions for all $i=1,\hdots,q$, Equation (\ref{cova_mult}) defines a family of matrix-valued isotropic covariance functions. Although this is a special case, it provides a useful multivariate extension of a subclass of the isotropic covariance models introduced by \cite{anderes2020isotropic}. It also complements previous results on vector-valued random fields on linear networks by \cite{filosi2025vector}, where isotropic matrix-valued covariance functions on trees were derived.

\section{Weighted local likelihood}
\label{sec:simulation}

\subsection{Description of the method}

We consider the weighted local likelihood estimation proposed by \cite{anderes2011local}. We now describe the method and the necessary modifications required to adapt it to the network context. 

Let $X(s_1)\hdots,X(s_N)$ be a realization of a GRF on $\mathcal{L}$. For any given site $s\in\mathcal{L}$, let $\mathcal{N}_{s,k}$ denote the set of $k$ observations nearest to $s$ with respect to the resistance metric. Recall that our goal is to estimate $a$ and $b$ at the selected site $s$. Given weights $\lambda_1,\hdots,\lambda_N$, the weighted local likelihood function is defined as
\begin{equation*}
\label{eq:local_lik}
\mathcal{W}\big(a_0,b_0 \mid X(s_1),\hdots,X(s_N)\big) = \sum_{k=1}^N \lambda_k \big[  \text{Lik}(a_0,b_0 \mid \mathcal{N}_{s,k}) -  \text{Lik}(a_0,b_0 \mid \mathcal{N}_{s,k-1}) \big],
\end{equation*}
where $\text{Lik}(a_0,b_0 \mid \mathcal{N}_{s,k})$ denotes the log-likelihood of the data in $\mathcal{N}_{s,k}$, with $a_0$ and $b_0$ representing fixed values of the functions $a(\cdot)$ and $b(\cdot)$. As in \cite{anderes2011local}, we define $\text{Lik}(a_0,b_0 \mid \mathcal{N}_{s,0})$ as 0. 

The local estimate of $\left(a(s),b(s)\right)$ is given by
$$  \big(\hat{a}(s),\hat{b}(s)\big) = \text{arg} \max_{(a_0,b_0)}\mathcal{W}\big(a_0,b_0 \mid X(s_1),\hdots,X(s_N)\big).   $$
The principle behind this method is that we locally fit an isotropic model, not simply by taking a window of data around $s$ (a strategy also used in \citealp{anderes2011local}, which has been shown to be less competitive), but by considering inflating neighborhoods and gradually assigning less weight to those that cover a very large region. For this reason, the weights include a bandwidth parameter to control this effect.

\subsection{Simulation}

In our simulation experiments, we focus on two main aspects. First, we study the performance of the weighted local likelihood method in the linear network setting, where it has not been previously examined. Second, we evaluate both the statistical and computational performance of compactly supported weights $\lambda_k^{\text{cs}} $ in comparison with Gaussian weights $\lambda_k^{\text{g}} $. 

Specifically, we consider the following weighting schemes:
\begin{equation}
\label{weight:cs}
\lambda_k^{\text{cs}} = \left(  1 - \frac{\text{d}_R(s,s_k)}{\tau}\right)_{+}, \qquad k=1,\hdots,N,
\end{equation}
where $(r)_+ = \max\{0,r\}$ denotes the positive part function, and 
\begin{equation}
\label{weight:gauss}
\lambda_k^{\text{g}} = \exp\left(-\frac{\text{d}_R(s,s_k)^2}{\eta^2}\right),  \qquad k=1,\hdots,N,
\end{equation}
where we assume that the sites $s_1,\hdots,s_N$ are ordered by their distance to $s$, for simplicity. Here, $\tau$  and $\eta$ are bandwidth parameters. Because compactly supported weights vanish beyond a certain distance when the bandwidth is moderate, only a subset of the neighborhoods $\mathcal{N}_{s,1},\hdots,\mathcal{N}_{s,N}$ may contribute to the local likelihood at $s$, which can reduce the computational cost.

\begin{figure}
    \centering
        \includegraphics[scale=0.8]{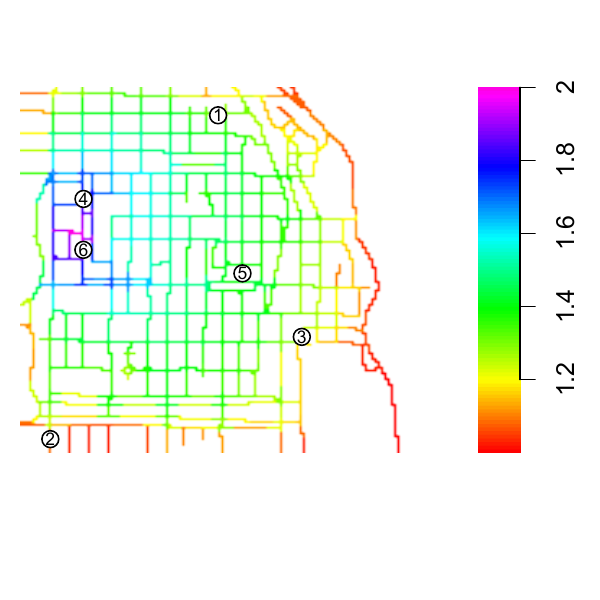}
    \includegraphics[scale=0.65]{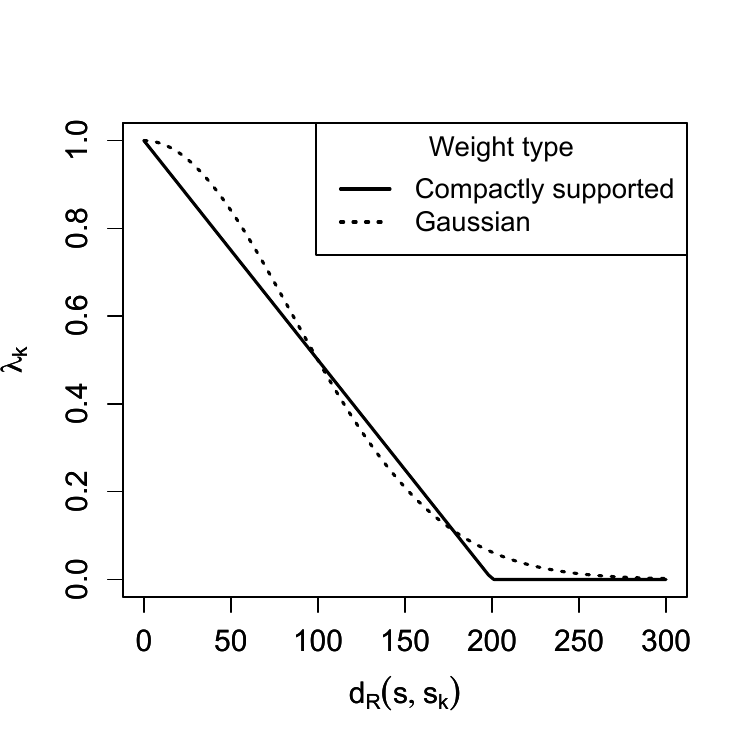}
    \caption{(Left) The function $a$ used in our simulations on the \texttt{Chicago} street network. Six randomly selected locations where $a$ will be estimated are highlighted. (Right) Decay of the weights as a function of the resistance metric for the compactly supported ($\tau=120$) and Gaussian ($\eta=120$) weighting schemes.}
    \label{fig:aa_sim}
\end{figure}

We consider the \texttt{Chicago} street network available in the \texttt{R} package \texttt{spatstat}, and the covariance function in Equation (\ref{cauchy}), with the parameterization discussed in Section \ref{sec:illustration}. Following \cite{anderes2011local}, we assume for simplicity that the variance of the random field is known and equal to one, so that we focus on the estimation of the local practical range. The function $a$ used in our study is shown in the left panel of Figure \ref{fig:aa_sim}. We also highlight six randomly selected locations on the network at which $a$ will be estimated. 

 We simulate 100 independent realizations of the GRF on this domain at 503 locations (one at the midpoint of each edge). For each realization, we estimate $a$ at the six selected sites. The weighted local likelihood method is implemented with the weighting schemes described above. We use bandwidth parameters $\eta=120$ and $\tau=200$, which produce weight functions with a comparable spatial extent but distinct structural features (see the right panel of Figure \ref{fig:aa_sim}). For reference, the maximum resistance distance between two points in this network is approximately 675.

\begin{figure}
    \centering
    \includegraphics[scale=0.8]{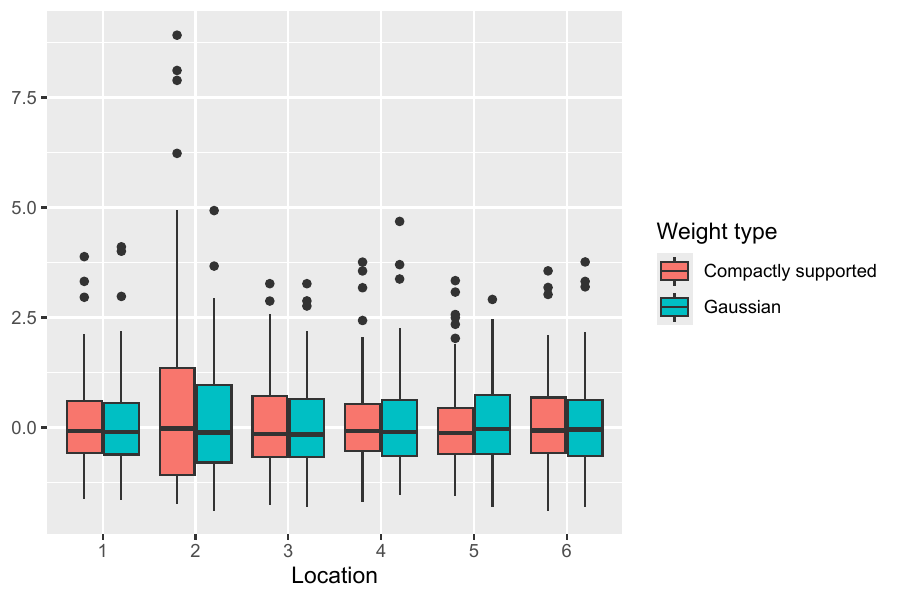}
    \caption{Centered boxplots of the estimates of $a$, at the six selected sites, across 100 independent realizations of a GRF with covariance function (\ref{cauchy}). Results are shown for the compactly supported and Gaussian weighting schemes.}
    \label{fig:boxplots}
\end{figure}

In Figure \ref{fig:boxplots}, we show the centered boxplots of the estimates across the 100 realizations, for each weighting scheme and each of the six selected sites. Overall, no significant bias is observed in any of the boxplots, demonstrating the effectiveness of the method. One point of interest is the higher variance in the estimates at site 2. This effect is consistent with the site's proximity to the boundary of the network, which reduces the number of neighbors that effectively contributes to the local estimate. A similar observation was reported by \cite{anderes2011local}.

Next, we comment on the comparison between the two types of weights. We observe similar results for both approaches, with a slight advantage for the compactly supported weights in terms of variance at sites 4, 5, and 6. However, at the boundary site 2, the variance is more pronounced when using compactly supported weights. Overall, compactly supported weights appear to be competitive from a statistical perspective.

We also compare the computational performance of the weighting schemes by measuring evaluation times as a function of the sample size of the GRF. Specifically, we simulate GRFs with $N=503, 1006, 1509$, and 2012, which correspond to taking 1, 2, 3, and 4 locations per edge, respectively, and implement the method using compactly supported weights with bandwidth parameters $\tau=100,150$, and 200. For each configuration, we report the time required for a single evaluation of the weighted local likelihood function. 

The results are summarized in Table \ref{tab:time}, allowing us to assess the computational gains relative to the Gaussian weighting scheme. For $\tau=100$ and 150, the reduction in computation time is particularly pronounced. For instance, when $\tau=150$, the computation time is reduced by approximately $50\%$ compared with Gaussian weighting scheme, whereas for $\tau=100$, evaluating the objective function takes at most $2\%$ of the time required with Gaussian weights, and for most sample sizes the percentage is even lower.

\begin{table}
\caption{Times (in seconds) for a single evaluation of the weighted local likelihood function under each weighting scheme, as a function of the sample size $N$. For the compactly supported weights we consider $\tau=100,150$, and 200. The simulation was performed on an Apple M1 chip with 8GB of memory. }
\centering
\begin{tabular}{|c|c|c|c|c|} \hline 
                                     & $N=503$ & $N=1006$  & $N=1509$ & $N=2012$  \\ \hline 
Compactly supported ($\tau=100$)   &  0.046      & 0.322          & 1.348         & 5.037    \\  \hline 
Compactly supported ($\tau=150$)   &  1.056      & 20.39          & 107.8         & 372.7    \\  \hline 
Compactly supported ($\tau=200$)   &  2.009      & 35.41          & 196.7        & 648.6    \\  \hline 
  Gaussian                                         & 2.257       & 41.35          &  225.8        & 759.0   \\  \hline 
\end{tabular}
\label{tab:time}
\end{table}

Our findings both support the observations reported in \cite{anderes2011local} regarding reduced precision near the domain boundaries and provide additional insights into the statistical and computational performance of different weighting schemes. A promising strategy is to use a combination of weights, employing compactly supported weights in the interior and adapting the weights near the boundaries. We refer the reader to \cite{anderes2011local} for a more comprehensive treatment of weights and bandwidth selection.

\section{Discussion}
\label{sec:discussion}

This paper proposes non-stationary covariance functions for random fields defined on linear networks. The models exhibit local isotropy with respect to the resistance metric, allowing both the variance and correlation range to vary spatially. We provide a useful stochastic representation for these random fields and extend the covariance functions from the scalar to the matrix-valued setting, enabling the modeling of vector-valued random fields. We demonstrate the effectiveness of a weighted local likelihood approach with compactly supported weights for estimating the proposed covariance models, which provide a balance between computational efficiency and statistical accuracy.

Although constructive criticism has been raised regarding covariance functions based on the resistance metric \citep{bolin2024gaussian}, primarily because existing covariance families under this metric do not provide explicit control over the mean-square differentiability of the random field (that is, there is no direct analogue of the Mat\'ern family of covariance functions under the resistance metric), we argue that the models introduced in this work offer considerable potential for practical applications. In the context of point process modeling, and in particular for LGCPs, local flexibility in variance and correlation range is essential for capturing heterogeneous aggregation structures, even when the smoothness level of the latent GRF is kept fixed \citep{dvovrak2019quick}. The proposed covariance models are also applicable to geostatistical data, such as continuous measurements along stream networks \citep{ver2010moving}. Since they do not allow control over the smoothness of the random field, their applicability is restricted to phenomena exhibiting relatively rough short-scale variations across space. Extending the framework to include flexible smoothness represents a promising direction for future research.

Other strategies for modeling non-stationarity in Euclidean spaces, such as domain deformation or different variants of kernel convolutions of isotropic random fields, could potentially be adapted to linear networks (see \citealp{schmidt2020flexible} and references therein). Investigating these approaches represents a promising direction for future research. Another interesting avenue is the development of nonparametric and semiparametric estimation methods in the network setting, analogous to the approaches of \cite{chang2010semiparametric} and \cite{kidd2022bayesian} in Euclidean spaces.

\section*{Acknowledgement}

Alfredo Alegr\'ia is partially supported by the National Agency for Research and Development of Chile through grant ANID Fondecyt 1251154.

\section*{Appendices}
\appendix

\section{Proof of Theorem \ref{teo1}}
\label{app1}

A direct calculation shows that, for all $w > 0$, 
\begin{eqnarray}
\label{eq_aux0}
 \frac{1}{\sqrt{\alpha(s,t) + \text{d}_R(s,t)w}}  & = &  \frac{1}{\sqrt{\alpha(s,t)}} \int_{\mathbb{R}} \exp\left\{ - \frac{w \text{d}_R(s,t)  y^2}{2\alpha(s,t)} \right\}  \frac{1}{\sqrt{2\pi}} \exp\left\{ -\frac{y^2}{2}\right\}  \text{d}y, \nonumber\\
 & = &  \frac{1}{\sqrt{\alpha(s,t)}}  {\rm E}\left[ \exp\left\{ - \frac{w \text{d}_R(s,t) Y^2 }{2\alpha(s,t)} \right\}  \right],
\end{eqnarray}
 where the expectation in (\ref{eq_aux0}) is taken with respect to $Y \sim \mathcal{N}(0,1)$. Consequently, $K$ in Equation (\ref{cova}) can be written as
\begin{equation}
    \label{eq_aux}
    K(s,t) = \frac{\sqrt{b(s)b(t)}}{\sqrt{\alpha(s,t)}}  \int_0^\infty  {\rm E}\left[  \exp\left\{  - \frac{w \text{d}_R(s,t)  Y^2}{2\alpha(s,t)}  \right\} \right] \text{d}F(w), \qquad s,t\in\mathcal{L}.
\end{equation}

Recall that $\phi(\cdot; \mu, \sigma^2)$ denotes the probability density function of a normal distribution with mean $\mu$ and variance $\sigma^2$. Using a convolution argument, we obtain
\begin{equation}
\label{convolution}
   \int_{\mathbb{R}} \phi\left(z; \frac{A(s)}{\sqrt{2}}, \frac{a(s)}{4w} \right) \phi\left(z; \frac{A(t)}{\sqrt{2}}, \frac{a(t)}{4w}\right)  \text{d}z  = \frac{\sqrt{w}}{\sqrt{\pi \alpha(s,t)} } \exp\left\{ - \frac{ w \big[  A(s)-A(t) \big]^2}{2\alpha(s,t)}  \right\},
    \end{equation}
for every $w > 0$, with this expression conditioned on the values, $A(s)$ and $A(t)$, of the auxiliary GRF defined in Equation (\ref{auxiliary}). 

Notice that $A(s)-A(t) \sim \mathcal{N}(0, \mathrm{d}_R(s,t))$ and, consequently, $\big[ A(s)-A(t) \big]/ \sqrt{\mathrm{d}_R(s,t)}$ is equal in distribution to $Y \sim \mathcal{N}(0,1)$. Taking the expectation in (\ref{convolution})  with respect to $(A(s),A(t))^\top$, we obtain
\begin{equation}
\label{expectation}
{\rm E}\left[  \int_{\mathbb{R}} \phi\left(z; \frac{A(s)}{\sqrt{2}}, \frac{a(s)}{4w} \right) \phi\left(z; \frac{A(t)}{\sqrt{2}}, \frac{a(t)}{4w}\right)  \text{d}z  \right]  
= \frac{\sqrt{w}}{\sqrt{\pi \alpha(s,t)} } \, {\rm E} \left[ \exp\left\{ - \frac{ w \text{d}_R(s,t) Y^2}{2\alpha(s,t)}  \right\} \right],
\end{equation}
where the right-hand side of (\ref{expectation}) simplifies to an expectation with respect to $Y\sim\mathcal{N}(0,1)$. 
Therefore, using Equation (\ref{expectation}) and Fubini's theorem, $K$ in (\ref{eq_aux}) can be written as
\begin{equation}
\label{eq4}
K(s,t) 
= \sqrt{\pi} \, {\rm E}\left[   
\int_0^\infty  \int_{\mathbb{R}} \sqrt{b(s)b(t)}  \, \phi\left(z; \frac{A(s)}{\sqrt{2}}, \frac{a(s)}{4w} \right) \phi\left(z; \frac{A(t)}{\sqrt{2}}, \frac{a(t)}{4w}\right)  \text{d}z \, w^{-1/2} \,  \text{d}F(w)     
\right].    
\end{equation}
We conclude that for any $N \in \mathbb{N}$, $r_1, \dots, r_N \in \mathbb{R}$, and $s_1, \dots, s_N \in \mathcal{L}$,
\begin{multline*} 
\sum_{i=1}^N \sum_{j=1}^N r_i r_j K(s_i,s_j)  =   \sum_{i=1}^N \sum_{j=1}^N r_i r_j \sqrt{\pi} \, {\rm E}\bigg[   
\int_0^\infty  \int_{\mathbb{R}} \sqrt{b(s_i)b(s_j)}  \, \phi\bigg(z; \frac{A(s_i)}{\sqrt{2}}, \frac{a(s_i)}{4w} \bigg)  \\ \phi\bigg(z; \frac{A(s_j)}{\sqrt{2}}, \frac{a(s_j)}{4w}\bigg)  \text{d}z \, w^{-1/2} \,  \text{d}F(w)   \bigg],     
\end{multline*}
that is,
$$\sum_{i=1}^N \sum_{j=1}^N r_i r_j K(s_i,s_j)  = \sqrt{\pi} \, {\rm E}\left[   
\int_0^\infty \int_{\mathbb{R}} \left\{ \sum_{i=1}^N r_i \sqrt{b(s_i)} \, \phi\left( z; \frac{A(s_i)}{\sqrt{2}}, \frac{a(s_i)}{4w} \right) \right\}^2    \text{d}z \, w^{-1/2}\, \text{d}F(w)   \right] \geq 0.$$
Therefore, $K$ is positive semidefinite.

\section{Proof of Proposition \ref{prop:stochastic}}
\label{app2}

On the one hand, using the independence of the involved random variables and the fact that ${\rm E}[\epsilon]=0$, it is straightforward to see that the random field has zero mean. On the other hand, by independence and noting that ${\rm E}\big[\epsilon^2\big] =1$, the covariance function is given by
\begin{eqnarray}
\label{eq3}
{\rm cov}\big[  X(s),X(t) \big]  & = & {\rm E}\big[\epsilon^2\big] {\rm E}\left[     \frac{\sqrt{b(s)b(t)}}{g(Z)}  \left( \frac{\pi}{W} \right)^{1/2}  \phi\left( Z; \frac{A(s)}{\sqrt{2}}, \frac{a(s)}{4W} \right)   \phi\left( Z; \frac{A(t)}{\sqrt{2}}, \frac{a(t)}{4W} \right)          \right],  \nonumber \\
& = &   {\rm E}\left[ \int_0^\infty \int_{\mathbb{R}} \frac{\sqrt{b(s)b(t)}}{g(z)}  \left( \frac{\pi}{w} \right)^{1/2}  \phi\left( z; \frac{A(s)}{\sqrt{2}}, \frac{a(s)}{4w} \right)     \phi\left( z; \frac{A(t)}{\sqrt{2}}, \frac{a(t)}{4w} \right)  g(z) \,\text{d}z \, \text{d}F(w)   \right],  \nonumber \\
& = &   {\rm E}\left[ \int_0^\infty \int_{\mathbb{R}}  \sqrt{b(s)b(t)}   \left( \frac{\pi}{w} \right)^{1/2} \phi\left( z; \frac{A(s)}{\sqrt{2}}, \frac{a(s)}{4w} \right)     \phi\left( z; \frac{A(t)}{\sqrt{2}}, \frac{a(t)}{4w} \right) \,\text{d}z \, \text{d}F(w)   \right], 
\end{eqnarray}
where the expectation in the last two equations is taken with respect to $(A(s),A(t))$. In Equation (\ref{eq3}), we identify the expression obtained in Equation (\ref{eq4}) in the proof of Theorem \ref{teo1}. Therefore, we recover the covariance function (\ref{cova}).

\section{Proof of Theorem \ref{teo2}}
\label{app3}

We first consider the simplified case where $\beta_{ij}=1$ for all $i,j=1,\hdots,q$. In this case, the function in Equation (\ref{cova_mult}) reduces to $\tilde{K}_{ij}$ defined by
$$  \tilde{K}_{ij}(s,t) =   \displaystyle \int_0^\infty \frac{\sqrt{b_i(s)b_j(t)}}{\sqrt{ \alpha_{ij}(s,t) + \text{d}_R(s,t)w}} \text{d}F(w).
 $$
Using the same arguments as in the proof of Theorem \ref{teo1}, we obtain
\begin{equation}
\label{eq4_multi}
\tilde{K}_{ij}(s,t) 
= \sqrt{\pi} \, {\rm E}\left[   
\int_0^\infty  \int_{\mathbb{R}} \sqrt{b_i(s)b_j(t)}  \, \phi\left(z; \frac{A(s)}{\sqrt{2}}, \frac{a_i(s)}{4w} \right) \phi\left(z; \frac{A(t)}{\sqrt{2}}, \frac{a_j(t)}{4w}\right)  \text{d}z \, w^{-1/2} \,  \text{d}F(w)     
\right].
\end{equation}
 Note that this is a multivariate version of Equation (\ref{eq4}).  Let $N \in \mathbb{N}$, $\bm{r}_1, \dots, \bm{r}_N \in \mathbb{R}^q$, and $s_1, \dots, s_N \in \mathcal{L}$, where $\bm{r}_k = (r_{k1},\hdots,r_{kq})^\top$. Then, 
\begin{multline*}
\sum_{k=1}^N\sum_{l=1}^N   \bm{r}_{k}^\top \tilde{\bm{K}}(s_k,s_l) \bm{r}_l = \sum_{k=1}^N\sum_{l=1}^N \sum_{i=1}^q \sum_{j=1}^q   r_{ki} r_{lj} \tilde{K}_{ij}(s_k,s_l)  \\ =   \sqrt{\pi} \, {\rm E}\left[    \int_0^\infty \int_{\mathbb{R}} \left\{ \sum_{k=1}^N \sum_{i=1}^q r_{ki} \sqrt{b_i(s_k)} \, \phi\left( z; \frac{A(s_k)}{\sqrt{2}}, \frac{a_i(s_k)}{4w} \right) \right\}^2    \text{d}z \, w^{-1/2}\, \text{d}F(w)   \right] \geq 0.
\end{multline*}
In conclusion, the $Nq\times Nq$ matrix
\begin{equation}
\label{matrix1}
 \begin{pmatrix}   \tilde{\bm{K}}(s_1,s_1) & \cdots & \tilde{\bm{K}}(s_1,s_N)  \\
\vdots & \ddots & \vdots \\
\tilde{\bm{K}}(s_N,s_1) & \cdots & \tilde{\bm{K}}(s_N,s_N)  
  \end{pmatrix}
  \end{equation}
  is positive semidefinite, as expected. 
  
  For the general case, where the parameters $\beta_{ij}$ are not necessarily identically equal to one, we observe that
\begin{equation}
\label{matrix2}
   \begin{pmatrix}   {\bm{K}}(s_1,s_1) & \cdots &{\bm{K}}(s_1,s_N)  \\
\vdots & \ddots & \vdots \\
{\bm{K}}(s_N,s_1) & \cdots & {\bm{K}}(s_N,s_N)   
  \end{pmatrix} = \left[\bm{1}_N\bm{1}_N^\top \otimes \bm{\beta}\right]\circ\begin{pmatrix}   \tilde{\bm{K}}(s_1,s_1) & \cdots & \tilde{\bm{K}}(s_1,s_N)  \\
\vdots & \ddots & \vdots \\
\tilde{\bm{K}}(s_N,s_1) & \cdots & \tilde{\bm{K}}(s_N,s_N)
  \end{pmatrix},
  \end{equation}
where $\circ$ and $\otimes$ denote the Hadamard and Kronecker products, respectively, and $\bm{1}_N$ is an $N$-dimensional vector of ones. 

Since $\bm{\beta}$, $\bm{1}_N\bm{1}_N^\top$, and the matrix in (\ref{matrix1}) are positive semidefinite, and positive semidefiniteness is preserved under both Kronecker and Hadamard products, it follows that (\ref{matrix2}) is positive semidefinite.

\bibliography{sample}
\bibliographystyle{apalike}

\end{document}